\documentclass[12pt,leqno]{amsart}

\usepackage{amssymb}
\usepackage{latexsym}

\newcommand\dspace{\lineskip=2pt\baselineskip=18pt\lineskiplimit=0pt}
\dspace

\newcommand\la{\lambda}
\newcommand\ga{\gamma}
\newcommand\tz{T(z, \bar z, u)}
\newcommand\ep{\epsilon}

\newcommand\ml{M_a^{k,l}}

\newcommand\dtz{D^2_z}

\newcommand\zz{(z,\bar z)}

\newcommand\pt{\tilde P}
\newcommand\ksi{\xi}
\enspace

\newcommand\nl{\newline}
\newcommand\fz{{\text{$F(z,\bar z,u)$}}}
\newcommand\fzs{{\text{$F^*(z^*,\bar z^*,u^*)$}}}

\newcommand\gjk{\gamma_{jk}}
\newcommand\glk{\gamma_{lk}}
\newcommand\hh{{\mathbb{H}}_{k,\mu}}

\newcommand\cen{{\text{${\mathbb C}^n$}}}
\newcommand\al{\alpha}
\newcommand\fss{f^{**}}
\newcommand\gss{g^{**}}

\numberwithin{equation}{section}
\theoremstyle{definition}

\newtheorem*{proo}{Proof}

\begin{document}

\quad
\vskip 3.3cm

\title[Generalized models and local invariants of Kohn-Nirenberg domains
]{ Generalized models and
local invariants of \\[2mm] Kohn-Nirenberg domains
}
\author[M. Kol\'a\v r]{Martin Kol\'a\v r}
\address{\newline
 Department of Mathematical Analysis, Masaryk University\newline
Jan\'a\v ckovo n\'am.\ 2a, 602 00 Brno\newline
{\it E-mail}: {\tt mkolar@math.muni.cz}}
\thanks{The author was supported by a grant of the GA \v CR no.
201/05/2117.}

\begin{abstract}
The main obstruction for constructing holomorphic reproducing kernels
of Cauchy type on weakly pseudoconvex domains is the
Kohn-Nirenberg phenomenon, i.e., nonexistence of supporting functions
and local nonconvexifiability.
In this paper we give an explicit verifiable characterization 
 of weakly pseudoconvex but locally nonconvexifiable hypersurfaces
of finite type in dimension two. It is expressed
in terms of a
generalized model, which captures local geometry of the hypersurface
 both  in the
complex tangential and nontangential directions. As an application we obtain
 a new
class of nonconvexifiable pseudoconvex hypersurfaces with
convex  models.
\end{abstract}


\maketitle

\section{Introduction}

Although 
 local convexity  of the boundary of a weakly pseudoconvex domain
is not a biholomorphic invariant, it  is often used as
an assumption which provides useful  tools to study
biholomorphically invariant objects
 (see e.g.\ \cite{DF}, \cite{DM}, \cite{M1}). The
natural invariant condition for results of this kind is then just
that the domain be  locally biholomorphic to a convex domain.
We consider here the problem of  finding an explicit characterization of such
domains.

Let  $D \subseteq \mathbb C^n$ be a domain with  smooth boundary
$M$ and 
  $p$ be a point on  $M$. Let  $\psi$ be 
 a local defining function  in a neighbourhood $U$ of $p$.
Recall that  $D$ is
pseudoconvex in  $U$ if for all
 $q\in M \cap U$ the Levi form
$$\mathcal L_q(\zeta) = \sum_{i,j = 1}^n \frac{\partial^2 \psi}{\partial z_i \partial \bar z_j}(q)
 \zeta_i \bar
 \zeta_j  $$
is nonnegative
on the complex tangent space to $M $ at $q$, i.e. 
on complex vectors
 $\zeta = (\zeta_1, \zeta_2, \dots \zeta_n)$
satisfying
  $\sum_{i=1}^n
 \frac{\partial \psi}{\partial z_i}(q) \zeta_i  =
0. $
When $D$ is strictly pseudoconvex (the Levi form is positive definite
at each point),
 fundamental results  of Henkin and Ramirez provide
 a holomorphic reproducing  kernel, a direct
analog of Cauchy's integral kernel from one  variable.
The construction is based on  holomorphic supporting functions
$h_q(z)$,  defined in a neighbourhood
of each $q \in M$,   whose zero set intersects $\bar D$ only
at $q$.
The reciprocal of a supporting function is the
analog of  the
fundamental 
function $\frac1{z-q}$                   
from one complex variable. 
The existence of such functions is a consequence 
of the following stronger result, known as Narasimhan's lemma.

\vskip 2mm
\noindent{\bf Theorem 1.1} (Narasimhan){\bf .}{\em\  Let  $D \subseteq \cen$
be a smoothly bounded domain, which is strictly
 pseudoconvex in a
neighbourhood of a boundary point $p$. Then there exist
local holomorphic coordinates around $p$ in which $\partial D $
is strictly convex.}
\vskip 2mm
 
In 1973, 
 Kohn and Nirenberg 
  found  an example
which shows that local
 convexifiability  does not in general extend to
weakly pseudoconvex domains. It  disproved the popular  conjecture
 that pseudoconvexity is equivalent to
 local convexifiability.
Their 
example is the  domain
$$D_{0} =  \{ (z,w) \in \Bbb C^2 \ \vert \  Im\ w > | z |^8 + \frac{15}7 
|z|^2 Re
z^6\}$$
 which does not admit a holomorphic 
supporting  function at the origin, and therefore is not convexifiable.
As a consequence, in connection with the efforts to extend Henkin's construction 
to (at least some) weakly pseudoconvex domains, it became 
an important  problem to characterize  domains for which 
this phenomenon occurs. 

Local convexifiability  of finite type domains turns out to be
 much stronger 
 than mere existence of supporting functions in all neighbouring points.
The problem of explicit characterization of locally convexifiable domains
 in complex dimension two was considered  in \cite{K1} and \cite{K3},
 with  a satisfactory answer on the level of model domains. The
results are based
on  exact conditions computed for  analogs of the Kohn-Nirenberg
example (see 
Proposition C below). They  do not take into account
behavior of the defining function in the complex non-tangential
direction. Moreover,  convexifiability of the model domain is
necessary, but not sufficient for convexifiability of the domain
itself. Hence the  possibility of other  obstructions to
convexifiability remained open.

Our aim 
in this paper 
is to complete the results of \cite{K1}, \cite{K3},
by introducing a generalized model
domain which (up to a borderline case)  carries the information about local
convexifiability of the domain itself. In this way,  the problem which
apriori requires to consider the whole infinite dimensional group
 of local biholomorphisms is reduced to a simple finite
dimensional problem, which is  often solved by an application of
the explicit results of \cite{K1}.

As an application, we find a new class  of
 nonconvexifiable domains.
While all previous examples are on the level of model domains, we describe nonconvexifiable
pseudoconvex domains
 whose model domains at each point are convex.

In Section 2 we use weighted coordinates to introduce the generalized 
model domain.
The construction is similar  to Catlin's
definition of multitype in ${\mathbb  C}^n$
(see \cite{C}), where only complex tangential
variables are considered. The weight of the complex nontangential direction
 provides a local biholomorphic
invariant with rational values. 
Section 3 gives two complementary conditions characterizing local convexifiability.
The proof of the sufficient condition relies on weighted homogeneity of the 
generalized model. The main tool in proving the necessary condition 
is the technique used  in \cite{CM}  for analyzing 
 biholomorphic transformations in weighted coordinates.   
Section 4 reviews the explicit conditions for local convexifiabity 
derived in \cite{K1}, \cite{K3}. In Section 5 we describe 
a new class of examples of nonconvexifiable domains which have
 convex  models at each point.  

A part of this work  was done while the author was visiting
the ESI program  Complex Analysis,
Operator Theory, and Applications to Mathematical Physics.
He would like to thank the organizers,  Friedrich Haslinger,
Emil Straube and Hans Upmeier for their invitation and hospitality,
and for the support recieved from ESI.
  
\section{Generalized model domains}

We will consider a  pseudoconvex domain
$D\subseteq {\mathbb  C}^2$ with smooth
boundary $M$, and a point $ p \in M $ of finite type $k$, where
$k>2$. Recall that the type of a point measures the maximal
order of contact of $M$ with complex curves passing through $p$. Strongly pseudoconvex
points are
 points of type two.
 It follows from  pseudoconvexity that $k$ is an even integer. 

We will use
local holomorphic coordinates $(z, w)$, $ z = x+iy$, $w = u+iv$,
centered at $p$,
 such that the positive $v$-axis
 gives the inner normal direction to $D$ at $p$
and write the defining equation for $M$ as
$$ v = \fz\,.
$$
 It follows directly from the definition of finite type 
(see  \cite{FS})
that we can choose the coordinates
so that the  equation takes form
\begin{align}
v =&P_1(z,\bar z)+o(|z|^{k}, u)\,,
\\
\intertext{where $P_1$ is a real valued homogeneous polynomial of degree $k$}
P_1(z, \bar z)&= a_0|z|^k + 
  \sum_{j=2,\dots ,k}|z|^{k-j} Re(a_j z^j)
\end{align}
for some $a_j \in {\mathbb  C}$ and $a_0 \in {\mathbb  R}^+$.
The
hypersurface
\begin{equation}
M_D =  \{ (z,w) \in {\mathbb  C}^2  \mid v = P_1(z,\bar z)\}
\end{equation}
is called a model hypersurface to $M$ at $p$.
 It is determined uniquely up to a linear
change of variables and addition of  a harmonic term
$Re\; \al z^k$ (see e.g. \cite{K3}).

\vskip 2mm
\noindent{\bf Definition 2.1. } 
 $M$ is called locally convexifiable at $p $ if there
exist local holomorphic coordinates in a full neighbourhood
of $p$ such that $M$ is convex with respect to the underlying
linear space.
\vskip 2mm
The polynomial $P_1$ in (2.1) captures local behavior of $M$ in
the complex tangential direction. In order to study
convexifiability of $M$ near $p$ we wish to take into account also
 behavior in complex nontangential directions,  involving the
variable $u$. A natural tool for  this is provided by weighted
coordinates.
 We will assign weights to the coordinates
$z$, $w$ and $u$.

\vskip 2mm
\noindent{\bf Definition 2.2. } 
 A weight  vector
 is a pair of  rational numbers $\la = (\la_1,\la_2)$
such that
$\la_1 = \frac1n$ for some $n\in {\mathbb  N}$ and there exist  integers $k_1, k_2$, with $k_2 > 0$
such that $k_1\la_1 + k_2\la_2 = 1$.
\vskip 2mm
 The weight  of  monomials
$c_{ijk}z^{i}\bar z^j u^k$ and $d_{ij}z^{i}w^k$ is defined to be  $(i+j)\la_1 + k\la_2$ and
$i\la_1 + k\la_2$, respectively.
A real valued polynomial $Q(z,\bar z, u)$   is $\la$-homogeneous of weight $\gamma$
if it is a sum of
 monomials of weight $\gamma$, and similarly for a holomorphic polynomial $P(z,w)$.
The weight of a function $Q(z,w)$ is the lowest of the weights of the nonzero terms in
its Taylor expansion at the origin
(when working in weighted coordinates it is helpful to use the two dimensional lattice
$\mathbb Z^+ \oplus \mathbb Z^+ $, where  a 
 monomial $z^a\bar z^b u^c$ in the Taylor expansion of $F$  
is assigned the lattice point $(a+b, c)$).
   
Given  the  point $p \in M$,  we first assign weight to the variable $z$, equal to
 the reciprocal of the type of the point, 
$\la_1 = \frac1k$. The weight of $w $ and $u$ will be denoted by
$\mu$, and is determined as follows. We will call a weight vector  $\la =
(\frac1k , \la_2)$ allowable, if there exist local holomorphic
coordinates such that the defining equation can be written in the
form
\begin{align*}
v= P(z,\bar z, u) + o_{wt}(1)\,,\\
\intertext{where}
P(z,\bar z, u) = \sum_{ \frac{m}k + \la_2 l = 1}
\sum_{j=0}^m  a_{mjl}z^j\bar z^{m-j}u^l
\end{align*}
is  a $\la$-homogeneous polynomial of weight one,  and $o_{wt}(1)$
denotes terms of weight bigger than one with respect to
the weight vector $(\frac1k , \la_2)$.

Clearly, for any real $\delta >0$ there are only finitely many rational numbers $\la_2$ bigger than $\delta$
 such that
 $(\frac1k , \la_2)$ is a weight. When the set of allowable values for $\la_2$ is finite, we set
$\mu$ to be the smallest of them. If the set is not bounded away from zero,
we set $\mu = 0$.
Using the change of variables formula (3.3) below,  it is straightforward to determine 
the value of $\mu$.

When $\mu > 0$, we fix local holomorphic coordinates $(z,w)$
which correspond to the weight $\la = (\frac1k, \mu)$. Hence the
defining equation has form
\begin{equation}
v = \sum_{ \frac{m}k + \mu l = 1}
\sum_{j=0}^m  a_{mjl}z^j\bar z^{m-j}u^l + o_{wt}(1)\,,
\end{equation}
where, as before,  the leading term will be  denoted by $P(z, \bar z, u)$.
The hypersurface  given by
$$
v =  P(z,\bar z, u)
$$
will be called a generalized model  to $M$ at $p$.
If $\mu=0$, a generalized model is defined to be
$$
v =  P_1(z,\bar z)\, ,
$$
hence in this case it coincides with the standard model hypersurface.

\section{Characterizations of convexifiability}

 Let $\hh$ denote  the set of all real
valued $\lambda$-homogeneous polynomials in $z,u$ of weight one
 which are
harmonic in the $z$ variable. Clearly, $h(z,u) \in \hh$ if and
only if
$$
h(z,u) =Re( \sum_{ \frac{m}k + \mu l = 1}
\alpha_m z^m u^l)
$$
for some $\alpha_m \in {\mathbb  C}$.

For $(z,u) \in {\mathbb  C} \times {\mathbb  R} $ and $\zeta \in {\mathbb  C}$ we will use the following notation
for the value of the real  Hessian restricted to  the $z$-direction:
$$
D^2_z P(z,u;\zeta) =  \sum_{i,j = 1}^{2}
 \frac{\partial^2 P}{\partial x_i \partial  x_j}(z,u) \ksi_i
 \ksi_j\,,$$
where $z=x_1 + i x_2$ and $\zeta = \xi_1 + i \xi_2$. Since $P$ is
weighted homogeneous, positivity of $\dtz P $ is determined by
its restriction to the set $S_2 \times S_1$, where $ S_2= \{(z,u)
\in {\mathbb  C} \times {\mathbb  R} : \vert z \vert^2 + u^2 =1\}$ and $S_1=
\{\zeta \in {\mathbb  C}: \vert \zeta \vert =1\}$.

Our aim is to prove the following pair of conditions for convexifiability
of $M$.

\vskip 2mm
\noindent{\bf Proposition A. } {\em
If there exists $h \in \hh$ such that $D^2_z (P + h) > 0$
 on $S_2\times S_1$, then M is locally convexifiable.}

\vskip 2mm
\begin{proo}
 Let $h \in \hh $ be such a function and let $\tilde P
= P + h$. Considering first the condition on the u-axis, it
follows from the assumption  that $\frac1{\mu}(1-\frac2k) $ is an
even integer and that  $P$  contains a term $\vert z\vert^2
u^{\frac1{\mu}(1-\frac2k)}$ with a positive coefficient.  By
homogeneity, we have
\begin{equation}
 \dtz \pt(z,u;\zeta) \geq  \epsilon (|z|^{k-2} +
 \vert u \vert^{\frac1{\mu}(1-\frac2k)})\vert \zeta\vert^2
\end{equation}
for a sufficiently small $\epsilon > 0$.
Consider the change of
coordinates
\begin{align*}
w^*&=w+iw^2, \ \ \ \ \ z^*=z, \\
\intertext{and denote $v^* = \fzs$ the defining equation 
in the new coordinates. We get }
F^*(z,\bar z, u - 2u\fz) &= \fz
+u^2 - (\fz)^2\,. \\
\intertext{Comparing terms of weight $\leq 1$ on both sides we obtain}
\fzs &= \pt(z^*, \bar z^*, u^*)  + (u^*)^2 + o_{wt}(1)\,.
\end{align*}
We will show  that $F^*$ is a convex function in a neighbourhood of the origin.
Denote $P^*(z,\bar z, u) = \pt(z,\bar z, u) + u^2$
(dropping stars on variables)  and consider
its $3 \times 3$ Hessian matrix  with respect
 to the variables $x$, $y$, $u$. The
$2 \times 2$ submatrix of the Hessian formed by the derivatives
with respect to $x$ and $y$ satisfies (3.1), hence is positive semidefinite.
 It remains to prove that the determinant of the Hessian is also nonnegative in a neighbourhood of the origin.
We calculate the weights of the terms entering into the
determinant. On the one hand,  $wt(P^*_{uu}P^*_{xx} P^*_{yy}) =
wt(P^*_{uu}P^*_{xy} P^*_{yx}) = 2- \frac4k$. On the other hand,
$wt( (P^*_{xu})^2 P^*_{yy}) = 3 - 2\mu - \frac4k$ and the same for all
the remaining  terms. Since $k>2$ and $\mu < \frac12$, we have
$3-2 \mu - \frac4k > 2-\frac4k$. By (3.1), the determinant
is nonnegative in a neighbourhood of the origin.

\qed\end{proo}

\vskip 2mm
\noindent
{\bf Proposition B. } {\em If $M$ is locally convexifiable, then there exists
 $h \in \hh$
such that
$\dtz (P +h) \ge 0 $ on $S_2 \times S_1$.}

\begin{proo}
 We will assume
that  $\dtz(P+h)$ is not nonnegative on $S_2 \times S_1$ for any $h
\in \hh$, and show that $M$ cannot be locally convex in any other
coordinates. If $\mu = 0$, the claim follows from Lemma 3 in
\cite{K3}, so we assume $\mu > 0$. 
Consider a  holomorphic transformation
\begin{equation}
\begin{aligned}
  z^*&=z+ g(z,w)
\\
 w^*&=w+f(z,w)\,.
\end{aligned}
\end{equation}
We may restrict ourselves to transformations which  preserves
the above  description of $M$,  so  $v^*=0$ is tangent to $M$ at the origin and
the positive $v^*$-axis points inside.
Let $F^*$ be the function describing $M$  in new
coordinates. Substituting (3.2) into $v^*=\fzs$, we get the 
change of variables formula
relating
 $F^*$ and $F$, $f$, $g$,
\begin{equation}
\begin{aligned}  F^*(z+g(z,u+iF(z, \bar z, u)),&\; \overline{ z +
g(z,u+iF(z, \bar z, u))},u+Re\
f(z,u+iF(z, \bar z, u))  =\\ = F(z,\bar z, u) &+ Im\ f(z,u+Re \
f(z,u+iF(z, \bar z, u)).
\end{aligned}
\end{equation}
Without any loss of generality we  assume that the
linear part of the transformation is partly normalized and
together require that
$$f=0\,,\quad  g=0\,, \quad  f_z=0\,,\quad 
g_z=0\,\quad  f_w =0\quad 
\mbox{at}\quad z=w=0.
$$
Germs of transformations which satisfy this condition  form a group, which 
will be denoted by $\mathcal G$.

We will call  $g(z,w)$ in (3.2)
strictly  superhomogeneous
if $wt(g) >  wt(z) = \frac1k $. 
By $\mathcal G_1$ we will denote the subgroup of $\mathcal G$ formed
 by germs of transformations  (3.2)
for which $g$ is   strictly superhomogeneous.

Next  we show  that every $ G \in \mathcal G$  can be factored uniquely as
$G = \tau\circ S$, where $\tau \in \mathcal G_1$ and $S$ is a polynomial 
transformation of the form
\begin{equation}
\begin{aligned}
  z^{**}&=z^* +  \sum_{j=1}^{[\frac1{k\mu}]} \delta_j (w^{*})^j, \\
 w^{**}&=w^*\,. 
\end{aligned}
\end{equation}

Indeed, write $G \in \mathcal G$  as
\begin{align*}
z^*& = z + \sum_{j=1}^{\infty} \al_j w^j +
 \sum_{i=1,j=0}^{\infty}\beta_{ij}z^iw^j, \\
w^* &= w +  \sum_{j=2}^{\infty} \epsilon_j w^j +
\sum_{i=1,j=0}^{\infty} \gamma_{ij}z^i w^j\,
\end{align*}
and compose it with a transformation $S$ of the form (3.4).
We claim there exist uniquely determined numbers 
$\delta_j$ such that $G \circ S \in \mathcal G_1$. 
In order to see this, 
 denote the composition by
\begin{equation*}
z^{**}=z+ \gss(z,w), \ \ \ \ \  
 w^{**}=w+\fss(z,w),
\end{equation*}
 where
\begin{align*}
  g^{**}(z,w)  = \sum_{j=1}^{\infty} \al_j w^j +
 \sum_{i=1,j=0}^{\infty}\beta_{ij}z^iw^j +&
  \sum_{m=1}^{[\frac1{k\mu}]}\delta_{m}
 (w + \sum_{j=2}^{\infty} \epsilon_j w^j +  \sum_{i,j=0}^{\infty}
 \gamma_{ij}z^i w^j)^m
\\
f^{**}(z,w)&=  \sum_{j=1}^{\infty} \epsilon_j w^j +  \sum_{i=1,j=0}^{\infty} \gamma_{ij}z^i w^j \,.
\end{align*}

There is a unique way  to choose the coefficients $\delta_j$
so  that $g^{**}$ is strictly superhomogeneous, namely 
$\delta_1 = - \alpha_1$, $\delta_2 = - \al_2 - \delta_1 \epsilon_2$,
and so on up to index $\big[\frac1{k\mu}\big]$.
Since the inverse of $S$ is of the same form as $S$,
we obtain the claimed decomposition.

Now  consider the transformation by  a representant of a fixed element  $ G_0 \in \mathcal G$ and use the above decomposition,
$G_0 = \tau_0 \circ S_0$.
Let $\tau_0$ be  of the form (3.2) and 
consider (3.3) for terms of weight less
or equal to one.

Now we consider  two cases.
When 
$f$ contains only terms of weight greater or equal to one, 
then  $\tau_0$ preserves  form (2.4), merely adding  harmonic terms
to $P$.    So it's enough to consider the effect
of $S_0$, given by

\begin{equation}
\begin{aligned}
  z^{*}&=z +  \sum_{j=1}^{[\frac1{k\mu}]} \delta_j w^j, \\
 w^{*}&=w\,. 
\end{aligned}
\end{equation}

By assumption, there exists a smooth curve $\ga$ through the origin
in $\mathbb C \times \mathbb R $
 and a vector $\zeta_0 \neq 0$  such that in a neighbourhood of $p$ 
\begin{equation}
 \dtz P(z,u;\zeta_0) \leq  - \epsilon (|z|^{k-2} +
 \vert u \vert^{\frac1{\mu}(1-\frac2k)})\vert \zeta_0\vert^2
\end{equation}
for $(z,u) \in \ga$ and  a sufficiently small $\epsilon > 0$.
Consider terms of  degree less or equal to $s+2$ in $F$, where
$s = (1 - \frac2k ) \frac1{\mu}$.
They will be  denoted   by  $\tz$, so  $ F= \tz + o(s+2)$.
Here the ordinary (nonweighted) degree is used,  $o(s+2)$ denoting 
terms of degree greater than $s+2$.
$T$ contains all terms which after applying $S_0$ can influence
terms of weight less or equal to one in $F^*$.
The estimate (3.6) holds also for  $\dtz T$,
perhaps with a smaller $\epsilon$ and a  smaller neighbourhood of $p$.
 After applying $S_0$ (and dropping stars on variables) 
 we consider terms of weight less or equal to one in $F^*$,
 and denote them  by $Q^*$. 
By (3.3)
 we obtain, 
modulo terms 
$o_{wt}(1)$,
$$Q^*(z, \bar z, u) = T(z+\sum_{m=1}^{[\frac1{k\mu}]}\delta_{m}u^m,
 \bar z + \sum_{m=1}^{[\frac1{k\mu}]}\delta_{m}u^m, u).$$
This transformation only shifts lines parallel to the z-axis,  and has no effect on convexity 
in the z direction.
Hence, denoting by $\ga^*$ the image of $\ga$, we have

\begin{equation}
 \dtz Q^*(z,u;\zeta_0) \leq  - \ep' (|z|^{k-2} +
 \vert u \vert^{\frac1{\mu}(1-\frac2k)})\vert \zeta_0\vert^2
\end{equation}
on $\ga^*$
in a sufficiently small neighbourhood of the origin and  $\ep'>0$ sufficiently small.
By homogeneity, it follows that $F^*$ cannot be
locally convex.

Now consider the second case, when $f$
contains terms of weight strictly less then one.  Then the 
leading weighted homogeneous
 term of $F^{*}$, which we  denote $\hat P$,   is harmonic in $z$. 
Let $a_{m0l}z^mu^l$ be the nonzero monomial 
in this leading term
with the lowest power of $u$, which dominates 
$\dtz \hat P$ along the $z$ plane. Denote by $\nu = \frac{m}k + \mu l$ its weight,
and let $k'=\nu k$. 
By  harmonicity,  there is a  smooth curve $\ga$
in $\mathbb C \times \mathbb R$ and $\zeta_0 \neq 0$,
 such that in 
a neighbourhood of $p$
\begin{equation}
 \dtz \hat P(z,u;\zeta_0) \leq
 - \ep' (|z|^{k'-2} +
 \vert u \vert^{\frac1{\nu}(1-\frac2{k'})})\vert \zeta_0\vert^2
\end{equation}
for $(z,u) \in \ga$ and  a sufficiently small $\epsilon' > 0$.
After applying $S_0$,
we use the same argument as in the first case.
Consider terms of  degree less or equal to $s+2$ in $F$, where
$s = (1 - \frac2{k'} ) \frac1{\nu}$, 
and denote them  by  $T^*$. Again, 
the same estimate holds  for  $\dtz T^*$ as for $\dtz \hat P$.
 After applying $S_0$, 
 we consider terms of weight less or equal to $\nu$ in $F^{**}$,
 denoted by  $Q^{**}$. 
By (3.3)
 we obtain, 
modulo terms 
$o_{wt}(\nu)$
$$Q^{**}(z, \bar z, u) = T^*(z+\sum_{j=1}^{[\frac1{k\nu}]}\delta_{j}u^j,
 \bar z + \sum_{j=1}^{[\frac1{k\nu}]}\delta_{j}u^j, u).$$
On the image of $\ga$,  (3.8) holds for $Q^{**}$. 
By homogeneity, it follows that $F^{**}$ cannot be
locally convex.

\qed
\end{proo}

\section{Convexifiability and  basic invariants}

For applications of Propositions A and B we will need the explicit results
obtained in \cite{K1}.
For two even integers $k, l$ and a positive real number $a$  denote
\begin{align*}
M^{k,l}_a &=
\{(z,w)\in {\mathbb  C}^2\ \mid \ Im \ w = P^{k,l}_a(z)\}\,,\\
\intertext{where}
P^{k,l}_a \zz &= |z|^k +a|z|^{k-l}Re\, z^l\, .
\end{align*}
The Kohn-Nirenberg example is a hypersurface of this type, with $k=8$, $l=6$
and $a = \frac{15}7$.

\vskip 2mm
\noindent
{\bf Proposition C. }
{\em $\ml$ is convex
if and only if
$a \le \gamma_{lk}$, where
\begin{align*}
\glk &= \frac{k}{l^2-k}\\
\intertext{if $l^2\geq 3k-2$ and}
\glk &= \sqrt{\frac{(4k - l^2 - 4)k^2}{(4k-4)(k^2
-l^2)}}
\end{align*}
if $l^2\leq 3k-2$.
Moreover, if $l$ is not a divisor of $k$, then this condition
is equivalent to convexifiability of $\ml$.}

\medskip

Now we consider the general case, when the hypersurface is given by
\begin{equation} v=P_1(z,\bar z)+o(\vert z \vert^{k},u),\end{equation}
where $P_1(z, \bar z)$ is
a polynomial of the form (2.2).
For such domains we  define the basic invariants. 

\vskip 2mm
\noindent
{\bf Definition 4.1.}{
\ Let the model hypersurface to $M$ at $p$ be given by (2.3),
  and let $l$ be 
an even integer, $0<l<k$.  We will 
call the  real number
$\kappa_{l} = \frac{\vert a_l \vert}{a_0}$
the Kohn-Nirenberg invariant of order $l$.
}

\vskip 2mm
\noindent

The numbers  $\kappa_l$  
 can be viewed   as the 3-rd level local biholomorphic invariants,
following  the  signature of the 
Levi form and the  
type of the point. The constant $ \frac{15}{7}$ in the definition
of $D_0$ is an invariant of order six.

In  the general case we have a sufficient and a necessary condition for
local convexifiability (\cite{K1}, \cite{K3}).
\vskip 2mm
\noindent
{\bf Theorem 4.2. }{\em Let the model  at $p \in M$ be given by
(2.3).
If
$$\sum_{j=2,4,...,k-2}\gjk^{-1} \kappa_{j} < 1,$$
then $M$ is  locally convexifiable at $p$.
}

\vskip 2mm
\noindent

\vskip 2mm
\noindent
{\bf Theorem 4.3. }{\em Let the model  at $p \in M$
be given by  (2.3). If $M $ is locally convexifiable at  $p$, then
\nl
(i)$\  \kappa_{j} \leq\gjk $ for all $j>\frac{k}2$
\nl
and
\nl
(ii)$\ \kappa_{j} \leq 2 \gjk $ for all $j \leq
\frac{k}2$.
}

\vskip 2mm
\noindent

\section{Examples of nonconvexifiable domains with convex models}

Using Propositions  B and C we can  find pseudoconvex domains  which are not
convexifiable, although their model domains at each point are convex.

As an example, consider 
 a hypersurface $M$  given by $ v = \fz$, where  
$$ F(z, \bar z, u) = P_1(z, \bar z) + Q(z, \bar z)u^s + o_{wt}(1),$$
and  $Q$ is a nonzero homogenous polynomial of degree
 $m < k$ without harmonic terms (both $m$ and $s$ are even).
It is easy to verify that in this case $ \mu = \frac1s - \frac{m}{ks}$ and $
P(z, \bar z, u) = P_1(z, \bar z) + Q(z, \bar z ) u^s$
is homogeneous of weight one 
with respect to the weight vector $(\frac1k, \frac1s - \frac{m}{ks})$.
  
Since $D^2_z$ of $Q(z, \bar z) u^s $ dominates
 $D^2_z$ of $P$ along the u axis, we get from Proposition  B the
 following condition.
If $M$ is convexifiable, then there exist $\alpha, \beta \in \mathbb C$ such that 
$P_1(z, \bar z)  + Re \; \alpha z^k $ and $ Q(z, \bar z) + Re \; \beta z^m$ are convex
.
Let us remark that for verifying this  one can use the simpler  complex form of
the  convexity condition -
a function $f(z, \bar z )$ is convex if and only if 
$$| f_{zz} | \leq f_{z\bar z}.$$
More explicitely, 
for $ a > 0 $   denote
$P_a^{6,4}(z,\bar z) = |z|^ 6 + a |z|^2 Re\, z^4$
and consider the hypersurface $M_a$ defined by
$$ v = |z|^8 + P_{a}^{6,4}(z,\bar z)  u^2 + |z|^2 u^8\,.
$$
Using (3.3) we verify that   $\mu = \frac1k = \frac18.$
The generalized  model is given by $ v = P(z,\bar z,u)$, where
$$
P(z,\bar z,u) = |z|^8 + (|z|^ 6 + a |z|^ 2 Re z^4) u^2\,.
$$
To obtain the pseudoconvexity condition for $M_a$  in terms of $F$, we take
$r = F-v$ as the defining function and use the coordinate
expression of the Levi form,
\begin{align*}
\mathcal L&=\vert r_w\vert^2 r_{z\bar z}+ \vert r_z\vert^2 r_{w \bar w}
- 2Re(r_{z\bar w}r_{\bar z}r_w)\,,
\\
\intertext{which gives}
\mathcal L&=F_{z \bar z}(1+\vert F_u \vert^2)+
F_{uu}\vert F_z\vert^2 + 2Re\,  F_{zu}(F_u+i)F_{\bar z}\,.
\\
\intertext{ Hence}
\mathcal L &= \frac14 \Delta F + o_{wt}(1-\frac{2}{k})\,.
\end{align*}
It follows that $M_a$ is pseudoconvex if  $\Delta P_a^{6,4} > 0$
away from the origin,
which leads to  $ a < \frac95$.
In this case $ l^2 = 3k-2$ in Proposition   C, which implies that
 $M_a $ is not convexifiable if $ a > \frac35$.
Hence for $\frac35 < a < \frac95 $ the hypersurface $M_a$ is pseudoconvex and
 nonconvexifiable.
 On the other hand, the model at zero is the convex domain $v = |z|^ 8$,
and all other points in a neighbourhood of zero are strongly pseudoconvex,
hence the models are also convex.


\bigskip
\begin{thebibliography}{12}
\itemsep=2pt

\bibitem{C}
  Catlin, D.,
   {\em Boundary invariants of pseudoconvex domains},
Ann.\ Math.\ {\bf 120} (1984), 529--586.


 \bibitem{D}
 D'Angelo, J., {\em  Orders od contact, real
hypersurfaces and applications}, Ann.\  Math.\ {\bf  115}
(1982), 615--637.

\bibitem {CM} S.S.Chern and J.Moser: \textit{Real hypersurfaces in
complex manifolds},  Acta Math.  \textbf{133} (1974), p. 219-271

\bibitem{DF}  Diederich, K.\ and Forn{\ae}ss, J.\ E., {\em Support
functions for convex domains of finite type},
Math.\ Z.\ {\bf 230} (1999), 145--164.

\bibitem{DF2}
 Diederich, K.\ 
 Forn{\ae}ss, J.\ E.
{\em Lineally convex domains of finite type: holomorphic support functions
},  Manuscripta Math.
{\bf 112} (2003), 403--431.

\bibitem{DM}  Diederich, K.\  and  McNeal, J.\ D., {\em
 Pointwise nonisotropic support functions on
convex domains},   Progress  Math.\ {\bf 188} (2000), 184--192.

\bibitem{FS}  Fornaess, J.\ E.\   and  Stensones, B., {\em
Lectures on Counterexamples in Several Complex Variables},
 Princeton Univ.\ Press 1987.

\bibitem{K} Kohn, J.\ J., \textit{Boundary behaviour of
$\bar \partial$ on weakly pseudoconvex manifolds of dimension two},
 J.\ Differential  Geom.\  \textbf{6} (1972),  523--542.

\bibitem{KN} Kohn, J.\ J.\ and Nirenberg, L., {\em  A pseudoconvex
domain not admitting a holomorphic support function},
 Math.\ Ann.\  (1973), 265--268.

\bibitem{K1} Kol\' a\v r, M., {\em  Convexifiability and supporting
functions in ${\mathbb C}^2$}, Math.\ Res.\ Lett.\ {\bf 2} (1995),
505--513.

\bibitem{K3} Kol\' a\v r, M., {\em  Necessary conditions for
local convexifiability of pseudoconvex domains in ${\mathbb  C}^2$},
Rend.\ Circ.\ Mat.\ Palermo {\bf 69} (2002), 109--116.

\bibitem{K5} M.Kol\'a\v r : \textit{Normal forms for hypersurfaces of
finite type in
 $ \mathbb C^2$}, Math. Res. Lett. \textbf{12} (2005)
 p. 897-910

 \bibitem{M1}  McNeal, J.\ D., {\em  Estimates on the Bergman Kernels
on Convex Domains},  Adv.\ Math.\ {\bf  109} (1994),
108--139.

\bibitem{M2}
 McNeal, J.\ D., {\em  Uniform subelliptic estimates on scaled convex domains of
              finite type},
 Proc.\ Amer.\ Math.\ Soc.\ {\bf  130} (2002),
 39--47 (electronic).

\end{thebibliography}
\end{document}